\documentclass[10pt]{amsart}

\usepackage{amssymb, amscd, amsmath, amsthm,  graphicx, latexsym}

\def\zz{{\bf Z}}

\def\cals{\mathcal{S}}

\begin{document}
\thispagestyle{empty}
\title{Slice knots with distinct Ozsv\'ath-Szab\'o and Rasmussen Invariants}

\author{Charles Livingston}
\address{Indiana University, Bloomington, Indiana, 47405} 
\email{livingst@indiana.edu}

\thanks{Research supported by the NSF}
\maketitle
\thispagestyle{empty}
\begin{abstract} As proved by Hedden and Ording, there exist knots for which the  Ozsv\'ath-Szab\'o and Rasmussen smooth concordance invariants, $\tau$ and $s$, differ.  The Hedden-Ording examples have nontrivial Alexander polynomials and are not  topologically slice.  It is shown in this note  that a simple manipulation of the Hedden-Ording examples yields a topologically slice Alexander polynomial one knot  for which $\tau$ and $s$ differ.  Manolescu and Owens have previously found a concordance invariant that is independent of both $\tau$ and $s$ on knots of polynomial one, and as a consequence have shown that the smooth concordance group of topologically slice knots contains a  summand isomorphic to $\zz \oplus \zz$.  It thus follows quickly from the observation in this note that this concordance group contains a subgroup isomorphic to $\zz \oplus \zz \oplus \zz$.
\end{abstract}
\vskip.3in

The   Ozsv\'ath-Szab\'o and Rasmussen knot concordance invariants~\cite{os, ra}, $\tau$ and $s$,  are each powerful   invariants, sufficient to resolve the Milnor conjecture  on the 4-genus of torus knots.  It had been conjectured that   $\tau = s/2$, but recently Hedden and Ording~\cite{ho} provided a counterexample, showing that certain doubled knots, including $ D_+(T_{2,3}, 2)$ and $ D_+(T_{2,5}, 4)$, satisfy $\tau = 0$ and $s = 2.$   Here $D_+(T_{p,q}, t)$  denotes the $t$-twisted positive double of the $(p,q)$-torus knot.  These examples demonstrate the richness of these new invariants.  However, in and of themselves, they do not reveal any new structure of the concordance group.  For instance, although some of these examples are algebraically slice, all can be shown not to be even topologically slice  using Casson-Gordon invariants~\cite{cg}.  (See~\cite{g} for techniques that resolve the nonsliceness of these particular doubled knots, and~\cite{ki} for a general discussion of Casson-Gordon invariants and doubled knots.) 
 
 In this note it will be shown that the basic Hedden-Ording examples can be manipulated to yield a knot with Alexander polynomial one (and thus, by Freedman~\cite{fr}, a topologically slice knot) for which $\tau$ and $s/2$ differ.  
 
 The smooth concordance group contains a subgroup  $\cals$ consisting of topologically slice knots.  In~\cite{liv2} it was  shown that $\tau$ and $s/2$ agree and are nonzero on some polynomial one knots:  both invariants take value 1 on   $D_+(T_{2,3},0) $ and $ D_+(T_{2,5},0)$.  It followed that $\cals$ contains a $\zz$ summand.  Manolescu and Owens~\cite{mo} developed a concordance invariant $\delta$ and showed $\delta(D_+(T_{2,3},0)) \ne \delta(D_+(T_{2,5},0))$, and thus showed $\cals$ contains  a summand isomorphic to $\zz \oplus \zz$.  The example here shows the independence of $\tau$ and $s/2$ on $\cals$ and   it follows from a quick determinant calculation that $\tau, s/2$, and $\delta$, yield a summand isomorphic to  $\zz \oplus \zz \oplus \zz$.

  \section{An Alexander polynomial one knot with $\tau \ne s/2$.}

 The doubled knots $D_+(K,t)$ bounds  a Seifert surfaces built by adding two bands to a disk, one with framing $-1$ and the other with framing $t$.  With respect to the corresponding basis of the first homology of the surface, the Seifert matrix  is:
 
\[
\left(
\begin{array}{cc}
 -1 & 1   \\
 0 & t \\
\end{array}
\right).
\]
 
\noindent It follows   from~\cite{ho} that  the connected sum $$K =D_+(T_{2,3}, 2)\, \#\, D_+(T_{2,3}, 2)\,  \#\,  D_+(T_{2,5}, 4)$$ satisfies $\tau(K) = 0$ and $s(K) = 6$.  With respect to the natural Seifert surface and bases for homology, the Seifert form of $K$ is given by the matrix $V_1$ below.

\vskip.2in 
\centerline{$
V_1 = \left(
\begin{array}{cccccc}
 -1 & 1  &  0 & 0 & 0 & 0  \\
 0 & 2  &  0 & 0 & 0 & 0  \\
 0 & 0  &  -1 & 1 & 0 & 0  \\
0 & 0  &  0 & 2 & 0 & 0  \\
  0 & 0 &  0 & 0 & -1 &1 \\
    0 &0   &  0 & 0 & 0 & 4 \\
\end{array}
\right) \hskip.3in 
V_2 = \left(
\begin{array}{cccccc}
 3 & 0  &  1 & 0 & 0 & 4  \\
 -1 & 2  &  0 & 0 & 0 & 0  \\
 1 & 0  &  -1 & 1 & 0 & 0  \\
0 & 0  &  0 & 2 & 0 & 0  \\
  0 & 0 &  0 & 0 & -1 &1 \\
    4 &0   &  0 & 0 & 0 & 4 \\
\end{array}
\right)
$}

\vskip.2in

If the first band in the Seifert surface is cut and reattached after twisting and linking with the other bands, a knot  $J$ with the Seifert matrix $V_2$ above  can be constructed.

For a knot with Seifert matrix $V$, the Alexander polynomial $\Delta_K(t)$ is given by $\det( V - tV^t)$; using this, a calculation gives $\Delta_J(t) = t^3$, and thus $J$ has trivial Alexander polynomial.

According to Livingston and Naik~\cite{ln}, cutting a band in a Seifert surface for a knot  and reattaching it can change $\tau$ by at most $\pm 1$ and can change $s$ by at most $\pm 2$.  It follows that $\tau(J) \in \{-1, 0, 1\}$ and $s(J) \in \{4, 6\}$.   (Note: $s$ is even and is bounded by twice the genus).  Clearly, $\tau(J) \ne s(J)/2$.


\begin{thebibliography}{9999}
 
  
  
\bibitem[CG]{cg}  {\bf A.~Casson, C.~McA.~Gordon}, {\sl Cobordism of classical knots},
in A la recherche de la Topologie perdue, ed. by Guillou
and Marin, Progress in Mathematics, Volume 62, 1986. (Originally
published as an Orsay Preprint, 1975.)

  
 \bibitem[F]{fr} {\bf M.~Freedman}, {\sl The topology of four-dimensional manifolds},  J. Differential Geom.  17  (1982),   357--453. 
 
 
\bibitem[G]{g} {\bf P.~Gilmer},  {\sl Slice knots in $S\sp{3}$}, 
  Quart. J. Math. Oxford Ser. (2) 34  (1983),
305--322. 


 
\bibitem[HO]{ho}  {\bf M.~Hedden, P. Ording}, {\sl The Ozsv\'ath-Szab\'o and Rasmussen concordance invariants are not equal}, (2005), arxiv.org/math/0512348. 

\bibitem[K]{ki} {\bf S.-G.~Kim}, {\sl Polynomial splittings of Casson-Gordon invariants},  Math. Proc. Cambridge Philos. Soc.  138  (2005),   59--78.


 \bibitem[L1]{liv1}  {\bf C.~Livingston}, {\sl Splitting the concordance group of algebraically slice knots},   Geom. Topol. 7 (2003), 641-643. 
  
  
\bibitem[L2]{liv2} {\bf C.~Livingston}, {\sl Computations of the 
Ozsv\'{a}th--Szab\'{o} knot concordance invariant}, Geom. Topol.   8  (2004), 735--742.


\bibitem[LN]{ln} {\bf C.~Livingston, S.~Naik}  {\sl  Ozsv\'ath-Szab\'o and Rasmussen invariants of doubled knots}, (2005), arxiv.org/math/0505361.

\bibitem[MO]{mo} {\bf C.~Manolescu, B.~Owens}, {\sl A concordance invariant from the Floer homology of double branched covers},  arxiv.org/math.GT/0508065.

\bibitem[OS]{os} {\bf P. ~Ozsv\'{a}th}, {\bf  Z. ~Szab\'{o}}, {\it  Knot Floer homology and the four-ball
genus},  Geom. Topol.     7  (2003),  615--639.    
 
\bibitem[Ra]{ra} {\bf J.~A.~Rasmussen}, {\sl  Khovanov homology and the slice genus}, (2003), arxiv.org/math/0306378.
 
 
\end{thebibliography}
\end{document}